\documentstyle{article}
\newcommand{\R}{{{\rm I}\!{\rm R}}}
\newcommand{\C}{{{\rm C}\!\!\!{\rm I}\,}}
\renewcommand{\epsilon}{\varepsilon}
\renewcommand{\phi}{\varphi}

\newcommand{\eqd}{\buildrel {\rm def}\over =}
\newcommand{\mov}{\hspace{0.3cm}}
 
\begin{document}
 
\title{Zero sets of some classes of entire functions}
\author{{\sc Alexander Russakovskii}}
\date{}
{\parindent 0pt \maketitle}
 
\medskip
 
{\bf Abstract. } --  A method of constructing an entire function
with given zeros and estimates of growth is suggested. It gives
a possibility to describe zero sets of certain classes of entire
functions of one and several variables in terms of growth of
volume of these sets in certain polycylinders.
 
\medskip
 
{\bf R\'esum\'e. } -- On propose une m\'ethode pour la
construction de fonctions enti\`eres dont l'ensemble des zeros
est donne\'e avec des majorations de croissance. Cette m\'ethode
permet de d\'ecrire avec pr\'ecision pour certaines classes de
fonctions enti\`eres (d'une ou de plusieurs variables) la
croissance du volume des diviseurs de leurs z\'eros dans
certaines polycylindr\`es.
 
\bigskip
 
\begin{large}
Let $K$ be some set of entire functions in $\C^n, \mov n\geq 1.$
Denote by $D_f$ the divisor of zeros of a function $f(z)$ and by
$Z_K$ the set $\{D_f: \mov f\in K\}.$ We are interested in
descriptions of the sets $Z_K$ for certain classes $K.$
 
In this note we make a survey of a couple of results which
appeared in {\bf [Ru1], [Ru2], [RoRu]}. They are devoted mainly
to descriptions of zeros of entire functions $f(z)$ with the
property
 
\begin{equation}
M_f (h) \eqd sup \{|f(z)| : \mov |Im z | \leq h \} < \infty,
\mov \forall h > 0
\end{equation}
 
We denote the class of all entire functions satisfying (1) by
$B.$ This class is important, because it contains, for example,
entire characteristic functions of probability distributions in
$\R^n$ and Dirichlet series with imaginary frequencies.
 
Note, that even
in one variable the usual canonical products seem to be a
nonadequate tool for constructing functions with given zeros
satisfying (1), and one needs to develop special technique and
use appropriate characteristics of zero sets in this case.
 
Let $\Pi(r,h) \subset \C^n$ be a polycylinder (rectangular box
for $n=1$) $\{z: |Re z| \leq r, |Im z| \leq h \},$ and let $D$
be a divisor in $\C^n.$ Denote by $n_D (r, h)$ the volume of $D$
in $\Pi(r,h) $, i.e.
 
$$ n_D (r,h) = \int\limits_{D\cap \Pi(r,h) } (dd^c
|z|^2)^{n-1}.$$
 
For the case $n=1,$ the quantity $n_D(r,h)$ is just the number of
points of $D$ (counted with multiplicities) in $\Pi(r,h).$
 
A complete descriprtion of zeros (divisors) for the whole class
$B$ was obtained by I.P.Kamynin and I.V.Ostrovskii ({\bf [KO]}):
\begin{it}
a divisor $D$ in $\C$ belongs to $Z_B$ if and only if
\begin{equation}
\forall h > 0: \mov log n_D(r,h) = o ( r ), \mov r \to \infty
\end{equation}
\end{it}
 
Using this result, Kamynin and Ostrovskii give a complete
description of zero divisors of Hermitian - positive entire
functions (Fourier transforms of probability distributions in
$\R$). We denote this class by $H$ and remind that an entire
function $f(z), \mov z \in \C^n,$ is called Hermitian - positive
if $f(0) = 1$ and for all $x^{(1)} \in \R^n, \ldots, x^{(n)} \in
\R^n, \mov w_1 \in \C, \ldots, w_n \in \C,$
 
$$\sum\limits_{k,j=1}^{n} f( x^{(k)} - x^{(j)} ) w_k
\overline{w_j} \geq 0.$$
 
According to the result in {\bf [KO],}
\begin{it}
a divisor $D$ in $\C$ belongs to $Z_H$ if and only if the
following holds:
 
$(i) \mov D \cap \{ Re z = 0\} = \emptyset;$
 
$(ii) \mov D = - \overline{D}, \mov  where \mov
- \overline{D}= \{z: -\overline {z} \in D\};$
 
$(iii) D \mov satisfies \mov  (2).$
\end{it}
 
Our theorem below extends this result to the case of several
variables and hence both results together give an answer to a
question posed by Yu.V.Linnik and I.V.Ostrovskii in {\bf [LO]}.
 
\medskip
 
{\sc Theorem 1 }({\bf [RoRu]}).
\begin{it}
A divisor $D$ in $\C^n$ belongs to $Z_H$ if and only if
conditions (i) -- (iii) above hold.  \end{it}
 
\medskip
 
If we put restrictions on the growth of $M_f(h),$ the problem of
description of zero sets becomes more complicated. Let $\phi(t)$
be a nonnegative increasing convex $C^2$-function on $\R_+.$ By
$B_\phi$ we denote the class of functions $f(z)$ belonging to
$B$ and satisfying the estimate
 
$$log M_f(h) \leq C e^{C \phi (h)}$$
 
\noindent
for some $C=C(f).$
 
Note, that since $\phi$ is supposed to be convex, we are dealing
with functions of infinite order. We assume additionally, that
$\phi (2t) = O(\phi (t)).$
 
Denote by $h(t)$ the solution of the equation
 
$$h \phi (h) = t, \mov t \geq 0.$$
 
It is easy to see, that $h$ is an increasing function of $t$
with the properties
 
$$\lim_{t \to \infty} \frac{h(t)}{t} =\lim_{t \to \infty}
\frac{\phi(h(t))}{t} =0.$$
 
A description of divisors in $Z_{B_\phi}$ is given by the
following
 
\medskip
 
{\sc Theorem 2.}
\begin{it}
A divisor $D$ in $\C^n$ belongs to $Z_{B_\phi}$  if and only if
 
\begin{equation}
log n_D (r,h) \leq C (1 + \phi (h) + \phi (h (r)) )
\end{equation}
 
\noindent
for some $C > 0.$
\end{it}
 
If we compare (3) to (2), we see that $o(r)$ is replaced in (3)
by a more precise expression $\phi(h(r)).$
 
A typical example of $\phi$ satisfying our conditions is
$t^\rho, \mov \rho \geq 1.$ In this case $\phi(h(t)) =
t^{\frac{\rho}{\rho + 1}}.$ We formulate the corresponding
result below.
 
{\sc Corollary.}
\begin{it}
A divisor $D$ in $\C^n$ is a divisor of an entire function
$f(z)$ with
 
$$\limsup_{h \to\infty} {log log M_f(h) \over h^{\rho}} < \infty
\mov (\rho \geq 1)$$
 
if and only if for some $C>0$
 
$$log n_D(r,h) \leq C( 1 + h^\rho + r^{\frac{\rho}{\rho +
1}}).$$
 
\end{it}
 
For $n=1$ one can obtain more precise results by using a more
precise characteristic of the divisors. Namely, we allow our
boxes $\Pi(r,h)$ to move along the real axis and add one more
parameter setting
 
$$n_D (x;r,h) = n_{D-x} (r,h),$$
 
\noindent
where $D-x$ means a translation of $D.$
 
For simplicity we formulate one result for a subclass of
$B_{t^\rho}.$ We denote by $[\rho, \sigma]$ the class of entire
functions $f(z)$ of one variable satisfying the condition
 
$$ \limsup_{h \to\infty} \frac{log log M_f (h)}{h^\rho} \leq
\sigma,$$
 
\noindent
and by $[\rho, \sigma]^*$ a subclass in $[\rho, \sigma]$
consisting of functions possessing the property
 
$$\exists \epsilon \in (0,\rho), \mov \exists a,b > 0: \mov
\forall x \in \R,$$
 
$$\sup_{|x - t| \leq a} log |f(t)| \geq  -exp \left(b
|x|^{\frac{\rho - \epsilon}{1+ \rho - \epsilon}} \right).$$
 
\medskip
 
{\sc Theorem 3.}
\begin{it}
Let $D$ be a divisor in $\C$ and let $\rho > 1, \sigma \geq 0$
be given. Then $D$ belongs to $Z_{[\rho, \sigma]^*}$ if and only
if there exist $\epsilon \in (0,\rho), C>0$ and a function
$\sigma (t) \to \sigma \mov (t \to\infty),$ such that
 
$$log n_D(x;r,h) \leq C + r + \sigma(h) h^\rho + C
|x|^{\frac{\rho - \epsilon}{1+ \rho - \epsilon}}.$$
\end{it}
 
\medskip
 
{\sc Remark.} The condition $\phi (2t) = O(\phi (t))$ is not
necessary and is used only to simplify the formulations. Results
for more general scales of growth may be found in {\bf[Ru1]}.
 
\medskip
 
The main tool for the above results is provided by a theorem on
construction of entire function with the given divisor of zeros
and estimates of its growth in terms of $n_D(r,h)$ (or
$n_D(x;r,h)).$ Similar results for exhaustion of $\C^n$ by balls
or polydiscs are well-known (see, for example, {\bf [Lel],
[Lev], [LG], [Ro], [St]}). A remarkable result of H.Skoda (see
{\bf [Sk]}) may be used in the proofs of some (not all) of the
theorems above instead of our construction. In any case, it
seems that the method described below might be of interest
itself.  The general idea of the method was expressed by
L.I.Ronkin.  Actually, one goes along the lines of the classical
solution of Cousin's second problem, but with estimates. We
remind briefly Cousin's scheme.
 
Assume for simplicity that we have a divisor $D$ in a domain
$\Omega \in \C.$  Cover $\Omega$ by disks or squares $G_j.$ Take
a polynomial $P_j (z) = \prod_{w_k \in D\cap G_j}(z-w_k)$.  On
$G_{ij} \eqd G_i \cap G_j$ the ratio $P_i\over P_j$ is a
nonvanishing holomorphic function, hence having the form $exp
(g_{ij}).$ Suppose, we are able to find such holomorphic
functions $g_k,$ that $g_{ij} = g_i - g_j$ on $G_{ij}$ for all
$i,j.$ Then we can put $f(z) = P_j exp (- g_j)$ on $G_j,$ thus
obtaining a function with given zeros in $\Omega.$ However, to
be able to represent $g_{ij}$ in the form of the difference, we
need to have a cocycle condition
 
$$g_{ij} + g_{jk} + g_{ki} =0.$$
 
All we have is that
 
$$ \exp (g_{ij} + g_{jk} + g_{ki}) = 1,$$
 
\noindent
or
 
$$g_{ij} + g_{jk} + g_{ki} = 2\pi i \cdot N$$
 
\noindent
with $N = N_{ijk}$ integer.
 
The next step is to represent $N_{ijk}$ in the form
 
\begin{equation}
N_{ijk}  =  M_{ij} + M_{jk} + M_{ki}, \mov\mov M_{\alpha\beta}-
{\rm integers}, \end{equation}
 
\noindent
which is possible exactly when the second group of cohomologies
with integer coefficients in $\Omega$ is trivial, this being the
necessary and sufficient condition for the solvability of
Cousin's second problem. This condition holds for every domain
$\Omega$ in $\C,$ which is not true for an arbitrary domain in
$\C^n;$ however it obviously holds for $\Omega = \C^n.$
 
Once we have passed this step, we put $h_{ij} = g_{ij} - 2 \pi i
\cdot M_{ij}.$ For these fuctions the cocycle condition takes
place and, representing them in the form $g_i - g_j,$ we obtain
our solution.
 
The idea is to follow the described scheme with bounds on each
step.  We can have good bounds for polynomials $P_k$ in terms of
the number of points of the divisor in $G_k.$ Thus we get
estimates for $g_{ij}$ and $N_{ijk}.$ The crucial point is to be
able to solve the cohomological equation (4) with "good"
estimates of the solution.  Once this is reached, it remains
to give estimates for $g_i,$ which is done in a standard way
with the help of H\"ormander's $\overline{\partial} -$ methods
({\bf [H\"o]}).
 
It appears to be possible to give solution with "good" bounds
for (4).  Moreover, the solution may be chosen in such a way,
that the bounds include the number of zeros in elements of the
covering, situated along one direction (this is where the
quantities $n_D (x; r, h)$ appear).  This is verified by
evaluating the number of "free parameters" and assigning them
certain appropriate values recurrently. The precise statement
for one - dimensional situation may be formulated as follows.
 
\medskip
 
{\sc Theorem 4.}
\begin{it}
Let $D$ be a divisor in $\C$ and let
 
$$\log n_D (x; 1, |y|) \leq u(x + i y), \mov\mov \forall x \in
\R, \forall y \in \R,$$
 
\noindent
where $u(z)$ is a subharmonic function.
 
Then there exists such an entire function $f(z)$ with $D_f=D,$
that
 
$$\log \log  |f(z)| \leq C + 2\log (1+|z|^2) + \sup_{|w-z|\leq
1} u(w).$$
\end{it}
 
\medskip
 
In higher dimensions it is also possible to solve (4) with
bounds.  However, the number of "free parameters" allows us to
eliminate only one real direction, and we could obtain our
estimates in terms of, say, $n_D (x_1; r,h),$ which is enough
for the above results.
 
Some details in several variables must be modified. We are not
able to take polynomials any longer. However, appropriate "local
solutions" are M.Anderson's solutions of
$\partial\overline{\partial} -$ problem in the unit ball of
$\C^n$ (see {\bf [An]}). We only need to choose a center for the
ball far enough from the divisor. This can be done so that the
distance to the divisor is estimated from below in terms of its
volume, and we obtain the "right" estimates. The corresponding
statement which we don't formulate precisely looks similar to
theorem 4.
 
We conclude with a remark, that in one dimension the outlined
scheme may be used for constructing functions with given zeros
in domains of $\C$ with control of growth near the boundary of
the domain.
 
\bigskip
 
\begin{center}
{\bf References}
\end{center}
\end {large}
 
[An] \quad
M.Anderson. {\it Formulas for the $L^2 -$ minimal
solutions of the $\partial\overline{\partial} -$ equation in the
unit ball of $\C^n.$} Math. Scand., 1984.
 
[H\"o] \quad
L.H\"ormander.{\it Introduction to complex
analysis in several variables.} Van Nostrand, Princeton, 1966.
 
[KO] \quad
I.P.Kamynin, I.V.Ostrovskii. {\it On zero sets of
entire Hermitian - positive functions.} Siber. Math. J.,
{\bf 23}, 1982, p. 66 - 82.
 
[Lel] \quad
P.Lelong. {\it Fonctions enti\`eres (n variables)
et fonctions plurisousharmoniques d'ordre fini dans $\C^n.$} J.
d'Anal. Math., {\bf 12}, 1964, p. 365 - 406.
 
[LG] \quad
P.Lelong, L.Gruman. {\it Entire functions of several
complex variables.} Springer - Verlag, 1989.
 
[Lev] \quad
B.Ya.Levin. {\it Distribution of zeros of entire
functions.} GITTL, Moscow, 1956.
 
[LO] \quad
Yu.V.Linnik, I.V.Ostrovskii. {\it Decompositions of
random variables and vectors.} Nauka, Moscow, 1972.
 
[Ro] \quad
L.I.Ronkin. {\it Introduction to the theory of
entire functions of several complex variables.} Nauka, Moscow,
1971.
 
[RoRu] \quad
L.I.Ronkin, A.M.Russakovskii. {\it Zeros of entire
Hermitian - positive functions of several variables.} Analytical
Methods in Probability Theory and Operator Theory. Naukova
Dumka, Kiev, 1990, p. 25 - 33 (Russian).
 
[Ru1] \quad
A.M.Russakovskii. {\it Zeros of entire functions of
infinite order.} Theory of Functions, Functional Analysis and
Appl., {\bf 54}, 1990, p. 105 - 122 (Russian).
 
[Ru2] \quad
A.M.Russakovskii. {\it Description of zero sets of
a class of entire functions of several variables.} Operator
Theory and Subharmonic Functions. Naukova Dumka, Kiev, 1991, p.
121 - 125 (Russian).
 
[Sk] \quad
H.Skoda. {\it Solution a croissance du second
probleme de Cousin dans $\C^n.$} Ann. Inst. Fourier, {\bf 21},
1971, p. 11 - 23.
 
[St] \quad
W.Stoll. {\it Ganze Funktionen endlicher Ordnung mit
gegebenen Nullstellen.} Math. Z., {\bf 57}, 1953, p. 211 - 237.
 
\bigskip
 
\obeylines
Alexander Russakovskii
{\it Theory of Functions Department
Institute for Low Temperature Physics \& Engineering
47 Lenin Avenue
310164 Kharkov
Ukraine}
 
\end{document}